\documentclass{article}
\pdfoutput=1 
\usepackage{amssymb}
\usepackage{amsfonts}
\usepackage{amsmath}
\usepackage{amsthm}
\usepackage{latexsym}

\usepackage{mathtools}
\usepackage{tikz, float}
\usepackage[british]{babel}
\usepackage[hidelinks]{hyperref}
\usepackage{cleveref}
\usepackage{microtype}
\usepackage{algorithm}
\usepackage{algpseudocode}
\usepackage{enumerate}
\usepackage[section]{placeins}

\usetikzlibrary{decorations.pathmorphing}

\usepackage[all]{xy}

\newcommand{\rats}{\mbox{\(\mathbb Q\)}}
\newcommand{\reals}{\mbox{\(\mathbb R\)}}

\newtheorem{theorem}{THEOREM}

\newtheorem{proposition}[theorem]{PROPOSITION}

\newtheorem{definition}[theorem]{DEFINITION}

 \def\set#1{{\{ #1\}}}
 
 \def \tikzconst{2.4}
 
 \def\c#1{{\mathcal #1}}
 
\def\P{{\bf P}}
\def\F{{\bf F}}
\def\G{{\bf G}}
\def\H{{\bf H}}

\def\x{{\sf x}}
\def\y{{\sf y}}

\def\M{{\sf M}}
 
 \def\Cl{{\sf Cl}}
 
 \def\ll{{\sf l}}
 \def\bb{{\sf b}}
 \def\rr{{\sf r}}
 \def\tt{{\sf t}}
 
 \def\tbb{{triangle bi-boundary}}

\date{}
\title{Temporal logic of Minkowski spacetime}
\author{Robin Hirsch and Brett McLean}
 
\AtEndDocument{\bigskip{\footnotesize
\noindent  Department of Computer Science, University College London, Gower Street, London WC1E~6BT \par  
\addvspace{\smallskipamount}
\noindent  \textit{E-mail address}:  \texttt{r.hirsch@ucl.ac.uk} 

\addvspace{\medskipamount}
\noindent  Laboratoire J. A. Dieudonn\'e UMR CNRS 7351, Universit\'e Nice Sophia Antipolis, 06108~Nice Cedex~02 \par  
\noindent  \textit{E-mail address}: \texttt{brett.mclean@unice.fr}
}}
 
\begin{document}
\maketitle
 
\begin{abstract} 
\noindent We present the proof that the temporal logic of two-dimensional Minkowski spacetime is decidable, PSPACE-complete.    The proof is based on a type of  two-dimensional \emph{mosaic}. Then we present the modification of the proof so as to work for slower-than-light signals. Finally, a subframe of the slower-than-light Minkowski frame is used to prove the new result that the temporal logic of real intervals with \emph{during} as the accessibility relation is also PSPACE-complete.
\end{abstract}

\section*{Introduction}
The theory of relativity raises several interesting issues for temporal reasoning. To mention a few: in special relativity there can be no notion of absolute simultaneity and the separation of time from space adopted by conventional knowledge representation systems has to be rejected; in general relativity fundamental theorems of computation become false, for example the halting problem may be solvable when observers cross the event horizon of a black hole \cite{MH94}. 

     Logical, axiomatic treatments of relativity have been developed, in the main part based on first-order, or in some cases second-order,  logic; see \cite{AMN07} for references.    However, a fundamental aspect of relativity theory is the abandonment of an absolute spacetime frame, and the requirement that all measurements and observations should be taken relative to a localised viewpoint or frame of reference.    A modal treatment of relativity would therefore be more natural (in this respect) than a first-order one.    A second connection between relativity theory and modal logic is to consider the points of spacetime as the worlds of a Kripke frame with an accessibility defined by the possibility of sending a signal from one spacetime point to another.  This latter connection is the topic of this chapter.

According to relativity theory, it is impossible to send signals faster than the speed of light. However, the notion of a signal travelling at the speed of light is invariant under changes of reference frame. In contrast, travelling at `half the speed of light', say, is not invariant---the most precise frame-invariant statement is that such a signal is travelling slower than lightspeed. So we may consider the points of  spacetime as worlds in a Kripke frame where a spacetime point ${\sf x}$ is accessible from a spacetime point $\y$  if it is possible to send a signal from $\y$ to $\x$, equivalently, if $\x$ lies in the future light cone of $\y$. 

The normal  Minkowski spacetime considered in  special relativity has three spatial dimensions and one time dimension, and if we use light-seconds for units of distance,  a point $(x, y, z, t)$ is accessible from $(x', y', z', t')$ if $t\leq t'$ and $(x'-x)^2+(y'-y)^2+(z'-z)^2\leq (t'-t)^2$.  Goldblatt showed \cite{Gol80} that the \emph{modal} validities over this frame are axiomatised by $\mathbf{S4.2}$---reflexivity $\Box p\rightarrow p$,  transitivity $\Box p\rightarrow\Box\Box p$, and confluence  $\Diamond\Box p\rightarrow\Box\Diamond p$, along with instances of propositional tautologies, using modus ponens and generalisation as inference rules.    The  modal logic does not change if we vary the number of spatial dimensions, nor  does it change if we require that signals are restricted to travel at strictly less than the speed of light, indeed the class of arbitrary reflexive, transitive, confluent frames has the same modal logic.  The irreflexive version of the Minkowski spacetime frame has  also been axiomatised; the modal logic is $\mathbf{OI.2}$ \cite{SS02}. The \emph{exactly lightspeed} accessibility relation, however, yields an undecidable logic \cite{Shap10}.

The topic of this chapter is the \emph{temporal} logic of Minkowski spacetime: similar to the unimodal logic, but we use $\F$ instead of $\Diamond$ for `sometime in the future' and  we have a converse modality $\P$ to access spacetime points in the past.  With the extra expressive power of the temporal language,  the temporal validities of Minkowski spacetime now depend on the number of spatial dimensions \cite{HR18}.   We conjecture that the validities of Minkowski spacetime with at least three dimensions is not decidable, but this conjecture remains open.
The main case we consider here is the temporal logic of two-dimensional Minkowski spacetime, with a single space dimension and one time dimension.  A complete axiomatisation of the  temporal validities of 2D Minkowski spacetime is not known. However, we will see that there is an algorithm to determine whether a given temporal formula is valid over 2D Minkowski spacetime or not.

In two dimensions it simplifies things if we change our coordinate system by rotating the axes anticlockwise by $45^\circ$, so the point $(r, t)$ with spatial coordinate $r$ and temporal coordinate $t$ is represented as $(x, y)$ where $x=\frac 1 {\sqrt 2}(t+r), \; y=\frac 1{\sqrt 2}(t-r)$.  In this coordinate system a point $(x', y')$ is accessible from $(x, y)$ if and only if  $x\leq x'$ and $y\leq y'$. We write $(x, y)\leq (x', y')$ in this case.     We can view this frame as a kind of product frame.

Let $\c K_1=(W_1, R_1),\; \c K_2=(W_2, R_2)$ be two Kripke frames.  Then the  \emph{product frame} $(W_1\times W_2, R_1^*, R_2^*)$ has two accessibility relations, the first  for accessing points horizontally:  $(x, y)R_1^*(x', y')\iff(y=y'\wedge xR_1y)$, and the second vertically: $(x, y) R_2^*(x', y')\iff (x=x'\wedge yR_2y')$.  
  These products have been studied intensely \cite{Kurucz07,GKWZ03}.    The logic of the product frame ${(\reals, \leq)}\times{(\reals, \leq)}$ is undecidable \cite{RZ01}.
But in 2D Minkowski spacetime, observers cannot tell whether a light signal arrives from the left or the right, since the map $(x, y)\mapsto (y, x)$ is  an automorphism of $\reals^2$ equipped with the up-to-speed-of-light accessibility relation.  So here we consider a product frame $\c K_1\otimes\c K_2$ with worlds $W_1\times W_2$ and a single accessibility relation $R$ given by $(x, y)R(x', y')\iff (xR_1x'\wedge yR_2 y')$ (for reflexive relations, this is the composition of $R_1^*$ and $R_2^*$).
The 2D Minkowski frame we looked at above can be written as $(\reals, \leq)\otimes (\reals, \leq)$.  The product of irreflexive frames $(\reals, <)\otimes(\reals, <)$ is the frame of two-dimensional spacetime points under the accessibility `can send a signal to a distinct point at strictly less than the speed of light', which we will investigate later.  The restriction of this latter frame to $\set{(x, y): x+y>0}\subseteq\reals^2$  provides a frame for a temporal logic of the strict `during' relation between intervals with real endpoints \cite{SS02}. This will be the third and final frame we will investigate.

Let  $(\reals, \leq)\otimes(\reals, \leq)$ be the 2D frame of Minkowski spacetime we have described. In fact, we will work primarily with the frame $\M$ derived from $(\reals, \leq)\otimes(\reals, \leq)$ by making accessibility irreflexive (which is not the same as the frame $(\reals, <\nobreak)\otimes(\reals, <\nobreak)$).    An \emph{$\M$-model} consists of the frame $\M$ and a map from propositional letters to subsets of $\reals^2$.  A \emph{line} in $\M$ means a straight line. A line is \emph{light-like} if it has the form $x=c$ or $y=c$ for some constant $c$. A line is \emph{time-like} if is linearly ordered by accessibility. A line is \emph{space-like} if it has strictly negative slope, equivalently no   pair of distinct points from the line is ordered. For points $\x < \y \in \mathbb \M$, we write $[\x, \y]$ for the closed interval in the poset $\M$ bounded by $\x$ and $\y$. Geometrically, then, $[\x, \y]$ is a rectangle (modulo degenerate cases).

The problem we face, then, is to determine whether a given propositional temporal formula $\phi$ is valid over $\M$. Rather than trying to find sound and complete axioms for this logic (which remains an unsolved problem) we construct a decision procedure based on an attempt to construct a certain type of finite description of a model, which we call a boundary map.  The bulk of this work is based on \cite{HR18,HM18}.     These boundary maps can be thought of as a continuation of the research in \emph{mosaics} pioneered by N\'emeti and his group \cite[\S~4]{Nem95a}, and extended to the temporal logic of the reals with `until' and `since' in \cite{Rey11}.    

The idea is to start from  basic building blocks and to synthesise our models by combining them together in various ways.   The building blocks of N\'emeti's mosaics are typically one-dimensional, whereas the ones we use here look like two-dimensional rectangles.    We will use these rectangular blocks to produce an algorithm to determine validity of temporal formulas over $\M$. From this, an algorithm for  the reflexive frame $(\reals, \leq)\otimes(\reals, \leq)$ follows immediately, via the reduction given by applying $\P\psi \mapsto \psi \vee \P\psi$ and $\F\psi \mapsto \psi \vee \F\psi$ recursively to subformulas. Then, we will modify the algorithm to work for the frame $(\reals, <)\otimes(\reals, <)$.  A new extension presented here is to modify the procedure used for slower-than-light accessibility to work for the temporal logic of real intervals, where the accessibility relation is `during'.  In each case we can prove the decidability of the validity problem, and indeed these decision problems are all PSPACE-complete.

\section*{Temporal  filtration}
A propositional temporal formula is either a propositional letter $p, q, r, \ldots$ or is  built from smaller formulas using  either propositional connectives $\neg, \vee$ or unary temporal operators $\F, \P$.  We use standard  propositional abbreviations $\wedge, \rightarrow$ and temporal abbreviations $\G\phi=\neg\F\neg\phi,\;\H\phi=\neg\P\neg\phi$ (`always in the future'/`always in the past').

Given a temporal formula $\phi$, its closure set $\Cl(\phi)$ is the set of all subformulas and negated subformulas of $\phi$. The cardinality of $\Cl(\phi)$ is linear in $|\phi|$. 

Fixing a $\phi$ whose validity we want to decide, a \emph{consistent set} is a subset of $\Cl(\phi)$ whose formulas can be (simultaneously) satisfied at some point in some (arbitrary) temporal frame. A \emph{maximal consistent set} (abbreviated MCS) is a consistent set that is maximal with respect to set inclusion, i.e. for some point in a temporal frame it is \emph{all} the formulas of $\Cl(\phi)$ holding there. (So, in particular, an MCS includes   one but not both of a subformula and its negation, for each subformula of $\phi$.) We denote the set of all such MCSs by $\mathsf{MCS}$. 

It is well known that satisfiability in a temporal frame is $\mathsf{PSPACE}$-complete \cite{Spaan1993}. Hence it can be checked in polynomial space whether any given subset of $\Cl(\phi)$ is an MCS. Such a subset can be stored in a number of bits linear in $|\phi|$, by recording its indicator function.

 For $m, n \in \mathsf{MCS}$ we let 
\begin{align*}
m\lesssim n &\iff &&\forall\F\psi\in \Cl(\phi)\;  ((\psi\in n\rightarrow \F\psi\in m)\wedge(\F\psi\in n\rightarrow\F\psi\in m))    \\  
&&\wedge\,&
\forall\P\psi\in \Cl(\phi)\;   ((\psi\in m\rightarrow \P\psi\in n)\wedge(\P\psi\in m\rightarrow \P\psi\in n)).
\end{align*}
The Kripke frame  $(\mathsf{MCS}, \lesssim)$ will be finite and  transitive.    The set of MCSs on which $\lesssim$ is reflexive is partitioned into \emph{clusters}---maximal sets of MCSs in which every pair of  MCSs is  in $\lesssim$. We extend the notation $\lesssim$ to clusters, thus we may write $m\lesssim c$ for an MCS  $m$ and a cluster $c$ if $m\lesssim n$ for some (equivalently all) $n\in c$.  For $\psi\in \Cl(\phi)$ we may say that $\psi$ belongs to a cluster $c$ if $\psi\in m$ for some $m \in c$.  Observe that it is possible for both $p$ and $\neg p$ to belong to a cluster, though a cluster cannot have both a temporal formula (e.g. $\F\psi)$ and its negation.  %  The problem we face is to transform this canonical frame into an $\M$-model, so we wish to colour each point in $\M$ with an MCS in a truth-preserving way. 

A \emph{defect} of an MCS  $m$ is simply a formula $\F\psi\in m$ (a future defect) or a formula $\P \psi\in m$ (a past defect).   A formula $\F\psi$ is a defect of a cluster $c$ if it belongs to $c$ but $\psi$ does not belong to $c$ (with a similar definition for past defects).   If $m\lesssim n$, the formula $\F\psi$ is a defect of $m$, and either $\psi\in n$ or $\F\psi\in n$ then we say that the defect $\F\psi$ of $m$ is \emph{passed up} to $n$ (with a similar definition for clusters). 

Any $\M$-model defines a function $h:\reals^2\rightarrow \mathsf{MCS}$ recording the formulas in $\Cl(\phi)$ true at each point.  Any two distinct vertical lines and two distinct horizontal lines define a closed rectangle in $\reals^2$, and restricting $h$ to this rectangle we obtain an instance of what we call a \emph{closed rectangle model} (definition imminent).  Similarly we may define  open (or semi-open) rectangle models by omitting  all (or some) of the four boundary lines from the domain of $h$.

\begin{definition}[Rectangle model]\label{def:rect}
A \emph{rectangle model} $h:R\rightarrow \mathsf{MCS}$ has a rectangle  $R\subseteq\reals^2$ (with edges parallel to the coordinate axes) as its domain and satisfies
\begin{description}
\item[coherence]
$\x\leq\y\in R\rightarrow h(\x)\lesssim h(\y)$;
\item [no internal defects]
if $\F\psi\in h(\x)$ then either there exists $\y\in R$ with $\y>\x$ and $\psi\in h(\y)$ (no defect) or
\begin{itemize}
\item  $R$ includes the boundary point  $\y$ due east of $\x$ and $\F\psi\in h(\y)$ (defect passed east) or
\item  $R$ includes the boundary point $\y$ due north of $\x$ and $\F\psi\in h(\y)$ (defect passed north).
\end{itemize}
and (similarly) where all occurrences of $\P$ defects may be passed south or west.  
\end{description}
%   If for all $\x<\y\in R$ such that $h(\x)$ and $h(\y)$ inhabit the same cluster $c$, we have for all $m \in c$ there is $\z\in R$ with $\x<\z<\y$ and $h(\z)=m$, then we say that $h$ maps \emph{densely}.
\end{definition}
It is straightforward to prove a `truth lemma' stating that \emph{open} rectangle models are equivalent to $\M$-models on which a $\psi \in \Cl(\phi)$ holds at a point $\x$ precisely when $\psi \in h(\x)$.

 The plan is to provide a finite description of rectangle models by starting from simple rectangle models where a single cluster holds at all points in the interior of the rectangle.  Then we show how to combine these simple rectangle models in various ways to obtain more complex rectangle models, but still with  finite descriptions.  Finally we will obtain a finite description of an open rectangle model. If we manage this, with $\phi$ mentioned during the construction, then we have a finite description of an  $\M$-model for $\phi$.
 
 \medskip

Before continuing, it is worth noting that temporal formulas can describe quite complicated  $\M$-models.  Consider the purely modal formula (not using $\P$/$\H$)
\[ \phi=\F p_0\wedge\F p_1\wedge\G\bigwedge_{i\neq j\neq k<3}(((\F p_i\wedge\F p_j)\rightarrow\F p_k )\wedge 
(p_i\rightarrow \G\neg p_j)).\]
This formula is satisfiable at $(0, 0)$ in an $\M$-model with valuation $v$ if $v(p_0)$, $v(p_1)$, $v(p_2)$ are disjoint subsets of the line $x+y=1$ with a point in $v(p_i)$ between any pair of points in $v(p_j)$ and $v(p_k)$, where $(i, j, k)$ is any permutation of $(0, 1, 2)$. And all $\M$-models satisfying $\phi$  look roughly like that, i.e. there is a closed space-like line segment in the future where $v(p_0), v(p_1), v(p_2)$ cover disjoint non-empty subsets of the segment (they may also cover points off the spatial line), with the density property just described.    We can even  express that $v(p_0), v(p_1), v(p_2)$ covers the closed line segment (with no gaps) by the formula
\[ \psi=\phi\wedge \G(\F p\vee\P p \vee p)\]
where $p$ abbreviates $(p_0\vee p_1\vee p_2)$, though this property cannot be expressed by a purely modal formula.  In an $\M$-model of $\psi$ there is  a segment of a space-like line $l$  partitioned by $v(p_0), v(p_1), v(p_2)$ in which the convex subsets of $v(p_i)$ are closed subsegments of $l$ for $i<3$ and, by topological properties of $\reals$, uncountably many of these closed segments are singleton sets (single points on the line)---a topological deduction that works in $(\reals,\leq)\otimes(\reals, \leq)$  but not in $(\rats, \leq)\otimes(\rats, \leq)$. Based on this observation we can define  temporal formulas satisfiable in the latter frame but not the former \cite[Figure 14]{HR18}.

\medskip

% In a closed rectangle model we always have $\x\leq y\rightarrow h(\x)\lesssim h(\y)$, but also any defect $\F\psi$ of $h(\x)$ is either witnessed in the rectangle (there is $\y$ in the rectangle with $\x\leq\y$ and $\psi\in h(\y)$) or it is passed up to the MCS   holding at a point directly north or directly east of $\x$ on either the northern or eastern edge of the rectangle.  In other words, there are no internal defects in the rectangle, but defects may be passed up/down to other adjacent rectangles.  %Lets see how we can describe these rectangle models in a finite way.

We have seen that the behaviour of the MCSs occurring along a space-like line can be quite unruly, making them hard to describe finitely.  Time-like lines are much easier,  since the function $h$ is monotone with respect to $\leq$, so we may describe the MCSs holding at points along such a line by a \emph{trace}: a $\lesssim$-ordered  finite sequence  $(c_0, m_0, c_1, m_1, \ldots, m_{k-1}, c_k)$ where the $c_i$ are distinct clusters, and the $m_i$ are MCSs (illustrated in  \Cref{fig:trace}).   The clusters $c_0$ and $c_k$ are called the \emph{initial} and \emph{final} cluster of the trace, respectively.  If $t, t'$ are traces, $m\in \mathsf{MCS}$, where  $t\lesssim m\lesssim t'$ (i.e. each cluster/MCS  of $t$ is below $m$, which is below each cluster/MCS  of $t'$) we may join them together to form a single trace $t\oplus m\oplus t'$, by concatenation, with the proviso that if the final cluster of $t$ equals the initial cluster of $t'$ they are identified and $m$ is omitted.    A formula $\F\psi$ is a defect of a trace $t$ if either it is a defect of the final cluster of $t$ or it is a defect of some other MCS/cluster of $t$ not passed up to the following cluster/MCS  of $t$.

\begin{figure}[H]
\[
\xymatrix{\ar@{-}[r]^{c_0}&\bullet_{m_0}\ar@{-}[r]^{c_1}&\bullet_{m_1}\ar@{-}[r]&\ldots&\ar@{-}[r] &\bullet_{m_{k-1}}\ar@{-}[r]^{c_k}&}
\]
\caption{\label{fig:trace}A trace}
\end{figure}

All pertinent information about a cluster $c$ can be stored in a  number of bits linear in $|\phi|$, for we only need record the indicator functions of
\begin{itemize}
\item
a representative MCS $m \in c$ (for making $\lesssim$-comparisons), 
\item
$\{\psi \in \Cl(\phi) \mid \exists m' \in c : \psi \in m'\}$ (to know the defects of $c$).
\end{itemize}
 The maximal length of a chain of distinct clusters or irreflexive members of $\mathsf{MCS}$ is also linear in $|\phi|$ (because walking up such a chain, we permanently gain formulas of the form $\neg\F\psi$ or $\P\psi$ at a linear rate). Hence any trace can be stored using a quadratic number of bits, and the number of traces is exponential in $|\phi|$.

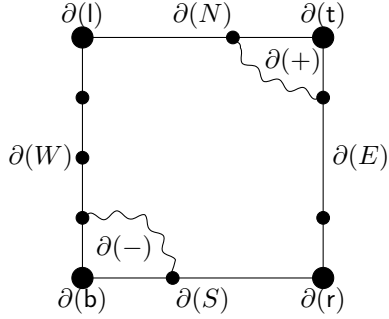
\begin{figure}[H]
\begin{center}

%\begin{tikzpicture}
\newcommand\Square[1]{+(-#1,-#1) rectangle +(#1,#1)}
\begin{tikzpicture}[scale=.8, decoration={snake, amplitude=1.5, segment length=10}]
\draw(0,0) \Square{2};

\node[below] at (-2,-2){$\partial(\bb)$};
\node [below] at (2,-2){$\partial(\rr)$};
\node [above] at (2,2){$\partial(\tt)$};
\node [above] at (-2,2){$\partial(\ll)$};

\node[left] at (-2,0) {$\partial(W)$};

\node[right] at(2,0){$\partial(E)$};
\node[above] at (0, 2){$\partial(N)$};
\node[below] at(0,-2){$\partial(S)$};

\draw[decorate](-2,-1).. controls  (-1,-.9) ..  (-0.5,-2);
\node at (-1.3,-1.5){$\partial(-)$};

\draw[decorate](2,1) .. controls (1, 1.3) .. (0.5,2);
\node at (1.5, 1.6){$\partial(+)$};

\draw [fill] (-2,-2) circle  [radius=5.0pt];
\draw [fill] (-2,2) circle [radius=5.0pt];
\draw [fill] (2,-2) circle [radius=5.0pt];
\draw [fill] (2,2) circle [radius=5.0pt];

\draw [fill] (-2,-1) circle [radius=3.0pt];
\draw [fill] (-2,0) circle [radius=3.0pt];
\draw [fill] (-2,1) circle [radius=3.0pt];

\draw [fill] (-0.5,-2) circle [radius=3.0pt];
\draw [fill] (0.5,2) circle [radius=3.0pt];
\draw [fill] (2,-1) circle [radius=3.0pt];
\draw [fill] (2,1) circle [radius=3.0pt];

\end{tikzpicture}
\end{center}

\caption{\label{fig:boundary}A boundary map}
\end{figure}

As indicated, the rectangles we will use have edges parallel to the coordinate axes, so their edges are segments of light-like lines, hence time-like lines (at least while we consider the case where signals may be sent at the speed of light).
Thus, we may describe the MCSs occurring along the perimeter of a rectangle $R$ in $\M$ by four  traces corresponding to the four open line segments bounding the rectangle and the four MCSs holding at the corners.   The distribution of MCSs in the interior of the rectangle can be complex. However, the only information we need to record is the minimal cluster holding at interior points arbitrarily close to the bottom corner and the maximal cluster holding at points arbitrarily near the top, as illustrated in \Cref{fig:boundary}.

\begin{definition}[Closed boundary map]\label{def:1}
A \emph{closed boundary map}  $\partial$  is a map from   $\set{N, S, E, W}\cup\set{\bb,\ll,\rr,\tt}\cup\set{-, +}$ to traces (for the first four), MCSs (for the next four), and clusters (for the last two), satisfying the following conditions.
\begin{itemize}
\item Temporal ordering must be respected, i.e. $\partial(\bb)\lesssim\partial(W)\lesssim\partial(\ll)\lesssim\partial(N)\lesssim\partial(\tt)$. And  (initial cluster of $\partial(W)) \lesssim \partial(-)\lesssim\partial(+)\lesssim$ (final cluster of $\partial(N)$).

\item Future defects of $\partial(+)$ must be passed up to either the final cluster of $\partial(N)$ or the final cluster of $\partial(E)$.  All future defects of $\partial(\bb)$ must be passed up to the initial cluster of $\partial(W)$ or the initial cluster of $\partial(S)$.

\item Dual conditions obtained from those above by  reflecting in either  or both diagonals, i.e. by swapping future/past, up/down,  $\tt/\bb,\; {\lesssim}/{\gtrsim}$,\/ N/W,\/ S/E, or by swapping N/E, S/W, $\ll/\rr$, or both, throughout.
\end{itemize}
\end{definition}

Some data can be omitted if some of the bounding edges and/or corners are not themselves included in the rectangle. 

\begin{definition}[Rectangular/rounded]\label{def:2}
Let  $\set{-, +}\subseteq D\subseteq\{-, +, N, S, E,\allowbreak W, \bb, \ll, \rr, \tt\}$.  If  a corner belongs to $D$ if and only if  the two adjacent edges are in $D$ (e.g. $\bb\in D\iff (S\in D\wedge W\in D)$) then we say that $D$ is \emph{rectangular}.   If we have only a one-way implication (whenever a corner is in $D$ then the two adjacent edges are also in $D$) we say that $D$ is \emph{rounded}.  
\end{definition}   The rectangular sets correspond to the different types of rectangles possible: open, partly open and closed rectangles, e.g. $\set{-, +}$ corresponds to an open rectangle, and $\{-, +, N, E,   \tt\}$ corresponds to a rectangle including its northern and eastern boundaries, but open to the south and west.    We may define general  boundary maps by restricting a closed boundary map to some rounded set $D$ by modifying \Cref{def:1}. So, for example, a future defect of $\partial(+)$ must be passed up to either the final cluster of $\partial(N)$ or of $\partial(E)$ (provided $N, E$ are in the set of $\partial$), but if neither $N$ nor $E$ is in the set then $\partial(+)$ can have no future defects.
%If $\tt\in D$ then   a future defect of $\partial$ is a future defect of $\partial(\tt)$ or a future defect of $\partial(N)$ or of $\partial(E)$.  When $\tt\not\in D$ then a future defect of $\partial$ means either a future defect of  $\partial(N)$ or of $\partial(E)$, or a future defect of $\partial(+)$ not passed up to the final cluster of either $\partial(N)$ or $\partial(E)$.   
%A past defect of $\partial$ is defined similarly.  Note, for an open rectangle where $D=\set{-, +}$, that a future defect is just a future defect of $\partial(+)$.    
Mostly we consider only rectangular sets, but when we eventually discuss interval logic we will need rounded sets including $S$ and $W$ but \emph{not} $\bb$.
%At one point later, we will also consider non-rectangular boundary maps $\partial$ called $\bb$-deleted boundary maps, where $W$ and $S$ are included in the domain of $\partial$ but not $\bb$.  A past defect of such a boundary map is a past defect of $\partial(W)$ or of $\partial(S)$.

As a notational convenience, we also define a \emph{one-point boundary map} to be a constant map $\partial : \{\bb, \tt\} \to \sf{MCS}$, which  may be identified with its image $m$. %When  $\partial(\tt)=\partial(\bb)$ the boundary map is called a \emph{one-point} boundary map and we may identify the boundary map by naming the MCS  $\partial(\bb)$.
Boundary maps that are not one-point boundary maps are said to be \emph{proper}.  

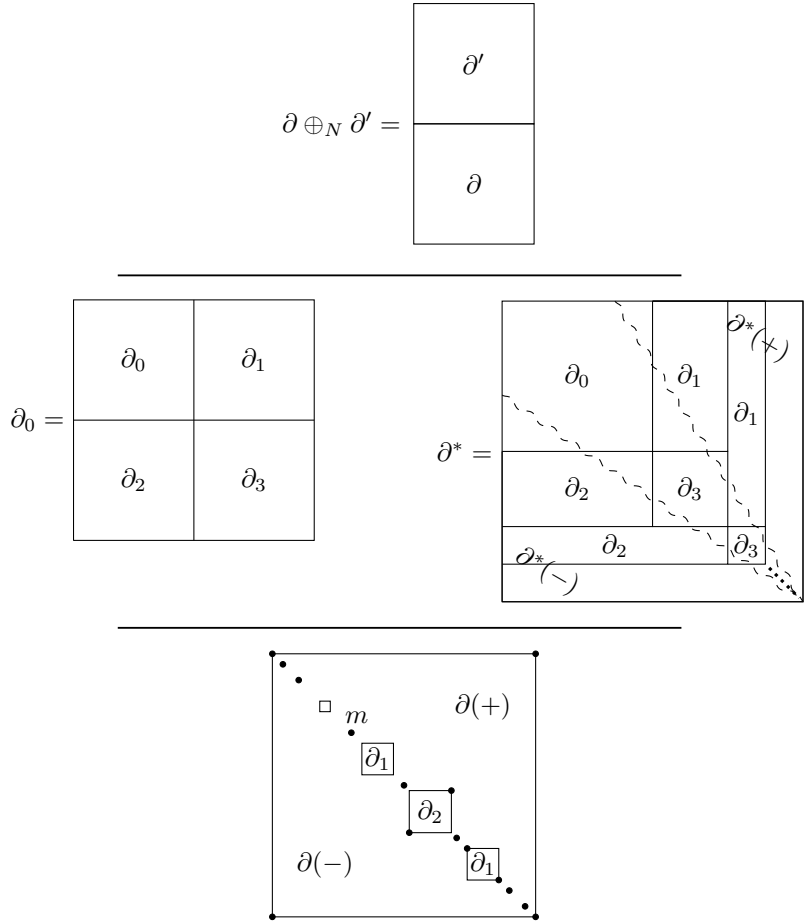
\begin{figure}  [H]

\begin{center}

\newcommand\Square[1]{+(-#1,-#1) rectangle +(#1,#1)}
\begin{tikzpicture}[scale=.4]

\draw(0,2) \Square{2}; 
\node at (0, 2) {$\partial'$};
\draw(0,-2) \Square{2};
\node at (0, -2){$\partial$};

\node at (-4.4,0) {$\partial\oplus_N\partial'=$};
\end{tikzpicture}\\
\rule{0.618\textwidth}{0.7pt}
\vspace{.1in}

\begin{tikzpicture}[scale=.4]

\draw(0, 0) --(0, 8)--(8,8)--(8, 0)--(0, 0);

\draw(0,4)--(8,4);
\draw(4,0)--(4, 8);

\node at (-1.2,4) {$\partial_0=$};
\node at (2,6){$\partial_0$};
\node at (6,6){$\partial_1$};
\node at(2, 2){$\partial_2$};
\node at (6,2){$\partial_3$};

\node at (0,-2.1){};

%\node at (.9, .3)      {\phantom   {$\partial_0(-)$}  };
%\node at (5, .3)  {\phantom  {$\partial_0(-)$}  };

\end{tikzpicture}
\hspace{.5in}
\begin{tikzpicture}[scale=.5, decoration={snake, amplitude=1}]
\node at (-1, 4){$\partial^*=$};
\draw(0, 0) --(0, 8)--(8,8)--(8, 0)--(0, 0);

\draw (0, 4)--(4,4)--(4, 8);

\node at (2, 6) {$\partial_0$};

\draw (0,2)--(4, 2)--(4, 4); \node at (2, 3) {$\partial_2$};
\draw (4,4)--(6, 4)--(6, 8)--(4, 8); \node at (5, 6) {$\partial_1$};
\draw (4, 2)--(6, 2)--(6, 4); \node at (5, 3) {$\partial_3$};

\draw (0, 1)--(6,1)--(6,2); \node at (3, 1.5) {$\partial_2$};
\draw (6,1)--(7,1)--(7, 8)--(6, 8);\draw(6, 2)--(7, 2);
\node at (6.5, 5) {$\partial_1$};
\node at (6.5, 1.5) {$\partial_3$};

\draw (0, 4)--(0,0)--(8, 0)--(8, 8)--(4, 8);

\draw[dashed, decorate](0,5.5)--(8,0);
\draw[dashed, decorate](3,8)--(8,0);
\node[rotate = -30] at (6.8, 7) {$\partial^*(+)$};
\node[rotate = -30] at (1.2, .8) {$\partial^*(-)$};

%\node[above] at (0, 8) {$\ll$};
%\node[below] at (8, 0){$\rr$};
%\node at (7.5, 0.5) {$\partial^*$};

\draw[dotted, very thick] (7.1,.9)--(7.9,0.1);

%\node at (.9, .3) {$\partial_0(-)$};
%\node at (5, .3) {$\partial_0(-)$};
%\node at (6, .2) {$\partial_0(-)$};
%\node[rotate=90] at(7.7, 7){$\partial_0(+)$};
%\node [rotate =90]at(7.7, 3){$\partial_0(+)$};
%\node[rotate = 90] at(7.8, 1.5){$\partial_0(+)$};

%\node at (4, -1) {(a)};
\end{tikzpicture}\\
\vspace{-.08in}
\rule{0.618\textwidth}{.7pt}
\vspace{.1in}

\begin{tikzpicture}[scale=.7]
\draw (0,0) -- (5,0) --(5,5) --(0,5) -- (0, 0);
%\node[above] at (2.5, 5) {$\partial(N)$};
%\node[below] at (2.5, 0){$\partial(S)$};
%\node[left] at (0, 2.5){$\partial(W)$};
%\node[right] at (5, 2.5){$\partial(E)$};
%\node[above, left] at (0, 5){$\partial(\l)$};
%\node[above, right] at (5, 5){$\partial(\tt)$};
%\node[below, left] at (0, 0){$\partial(\bb)$};
%\node[below, right] at (5, 0){$\partial(\rr)$};
\draw (3, 2) \Square{0.4};
\node at (3, 2) {$\partial_2$};
\draw (2, 3)\Square{0.3};
\node at (2, 3) {$\partial_1$};
\draw(1,4)\Square{0.1};
\draw(4, 1)\Square{0.3};
\node at (4, 1) {$\partial_1$};

\node at (1.6,3.8){$m$};

\draw [fill] (4.3,0.7) circle [radius=1.5pt];

\draw [fill] (3.7,1.3) circle [radius=1.5pt];

\draw [fill] (4.5,0.5) circle [radius=1.5pt];
\draw [fill] (.5,4.5) circle [radius=1.5pt];
\draw [fill] (2.5,2.5) circle [radius=1.5pt];
\draw [fill] (1.5,3.5) circle [radius=1.5pt];
\draw [fill] (3.5,1.5) circle [radius=1.5pt];
\draw [fill] (.2,4.8) circle [radius=1.5pt];
\draw [fill] (4.8,0.2) circle [radius=1.5pt];
%\node [below] at (4.5, .5) {$m_1$};

%\node [below] at (2.3, 2.5) {$m_0$};

%\node [below] at (.5, 4.5) {$m_1$};
%\node [below] at (2.6,1.6) {$\partial_2(\b)$};
%\node [above] at (3.4,2.3) {$\partial_2(\tt)$};f
\draw [fill] (3.4, 2.4) circle [radius=1.5pt];
\draw [fill] (2.6, 1.6) circle [radius=1.5pt];
\draw [fill] (0,0) circle [radius=1.5pt];
\draw [fill] (0,5) circle [radius=1.5pt];
\draw [fill] (5, 0) circle [radius=1.5pt];
\draw [fill] (5, 5) circle [radius=1.5pt];

\node at (1, 1) {$\partial(-)$};
\node at (4, 4) {$\partial(+)$};

\end{tikzpicture}

\end{center}
\caption{\label{fig:combine} A northern join, a southeastern limit, and a shuffle of $\partial_1, \partial_2, m$}
\end{figure}

In the simplest kind of proper boundary map we have $\partial(-)=\partial(+)$, corresponding to   rectangle models where a single cluster holds over the entire interior of the rectangle.  More complex boundary maps may be \emph{fabricated} from simpler ones,   by joins, limits, and shuffles (see \Cref{fig:combine}). 

 To first help with the intuition, given some propositional valuation, if $R$ is a rectangle in the plane, and if $R$ can be divided into two adjacent rectangles $R_1, R_2$ sharing a common edge, then the boundary map determined by $R$ will be the join  of the boundary map determined by $R_1$ and that determined by $R_2$ (in Figure~\ref{fig:combine} we show a \emph{northern} join).  
 
 If $\x, \y_i$ are all incomparable ($i<\omega$), if the boundary map $\partial_0$ determined by the rectangle $[\x\wedge\y_i, \x\vee\y_i]$ is constant (as $i$ varies), if the sequence $\y_i$ converges in a  southeasterly direction to $\y$, and if the maximal cluster of $[\x\wedge\y, \x\vee\y]$ extends all the way south to $\y$ and the minimal cluster all the way east to $\y$, then the boundary map determined by $[\x\wedge\y, \x\vee\y]$ will be a southeastern limit of $\partial_0$.  Intuitively, the region labelled with neither the maximal nor minimal cluster must be pinched at $\y$, as sketched in the second section of \Cref{fig:combine}.
 
  To picture a boundary map $\partial$  that is a shuffle of  closed boundary maps $\partial_0, \partial_1, \ldots$ (including at least one one-point boundary map $m$),  imagine a rectangle in the plane  with copies of rectangles for $\partial_0, \partial_1, \ldots$ disjointly but densely distributed in the manner of a Cantor set, along the diagonal from $\ll$ to $\rr$, with gaps along this diagonal filled by  the one-point boundary map $m$, with $\partial(-)$  covering the part of the rectangle below the diagonal and the copies of $\partial_0, \partial_1, \ldots$, and with $\partial(+)$ covering the area above---see page~\pageref{p:shuffle} for a bit more detail of this construction.
\begin{description}
\item[Joins]  If $\partial'$ fits to the north of $\partial$ (i.e. $\partial(N)=\partial'(S)$ is defined, \/ all or none of $\partial(W), \partial'(W), \partial(\ll),\partial'(\bb)$ are defined, and if defined  $\partial(\ll)=\partial'(\bb)$,\/ and similarly for $\partial(E), \partial'(E), \partial(\tt), \partial'(\rr)$),  then we may form the join $\partial\oplus_N\partial'$ defined by \begin{align*}(\partial\oplus_N\partial')(-) &= \partial(-) &(\partial\oplus_N\partial')(+) &= \partial'(+)\\
\shortintertext{and when the right-hand sides are defined}
(\partial\oplus_N\partial')(N) &= \partial'(N) &(\partial\oplus_N\partial')(S) &= \partial(S)\\
(\partial\oplus_N\partial')(W) &= \partial(W) \oplus \partial'(W) &(\partial\oplus_N\partial')(E) &= \partial(E) \oplus \partial'(E)\\
(\partial\oplus_N\partial')(\ll) &= \partial'(\ll) &(\partial\oplus_N\partial')(\tt) &= \partial'(\tt)\\
(\partial\oplus_N\partial')(\bb) &= \partial(\bb) &(\partial\oplus_N\partial')(\rr) &= \partial(\rr)
\end{align*}
(see the first part of \Cref{fig:combine}).  The join operator $\oplus_E$ is defined similarly.  

\item[Limits] If  
\begin{itemize}
\item  $\partial_0=(\partial_2\oplus_E\partial_3)\oplus_N(\partial_0\oplus_E\partial_1)$, 
\item  $\partial^*$ agrees with $\partial_0$ on $\set{-, +, \ll, W, N}$ (if defined),
\item
$\partial_1(E)$ is the trace consisting of the single cluster $\partial^*(+)$, and
$\partial_2(S)$ is the trace consisting of the single cluster $\partial^*(-)$, 
%\item   the domain of $\partial^*$ includes $\set{S, E, \rr}$,
\item $\partial^*(-) = \partial_0(-)$ (necessarily equal to $\partial_2(-)$), and $\partial^*(+) = \partial_0(+)$ (necessarily equal to $\partial_2(+)$),
\item if $\partial^*(S)$ is defined, every future defect is passed up to $\partial^*(-)$ or $\partial^*(\rr)$, and if $\partial^*(E)$ is defined, every past defect is passed down to $\partial^*(+)$ or $\partial^*(\rr)$, 
%\item %automatic from fouth condition and condition that $\partial^* is a boundary map
% every past defect of $\partial_0(-)$ is passed down to the initial cluster of either $\partial^*(W)$ or $\partial^*(S)$, and every future defect of $\partial_0(+)$ is passed up to the final cluster of either $\partial^*(N)$ or $\partial^*(E)$, 
 \end{itemize}
 then $\partial^*$ is a \emph{southeastern limit} of $\partial_0$ using $\partial_1, \partial_2, \partial_3$. (See the second part of \Cref{fig:combine}, where the wavy lines indicate the implied boundaries of the regions labelled with $\partial^*(-)$ and $\partial^*(+)$.)  Northeastern limits are defined similarly.

\item[Shuffles]   If $\partial_1, \partial_2, \ldots, \partial_{k}$ are closed boundary maps including at least one one-point boundary map, every future defect of $\partial(-)$ is passed up to $\partial_i(\bb)$ for some $i\leq k$, every past defect of $\partial(+)$ is passed down to $\partial_i(\tt)$ (some $i \leq k$), every past defect of each $\partial_i(\bb)$ is passed down to $\partial(-)$, and every future defect of each $\partial_i(\tt)$ is passed up to $\partial(+)$ (all $i$), then $\partial$ is a shuffle of $\partial_1, \ldots, \partial_{k}$ (see the last part of \Cref{fig:combine}).
\end{description}

\begin{definition}[Fabricated]\label{def:fabricated}
A boundary map $\partial$ is \emph{fabricated} if it occurs in a sequence of distinct boundary maps each of which is either \emph{simple} ($\partial(-)=\partial(+)$) or obtainable as the join, limit, or shuffle of earlier boundary maps in the sequence.    
\end{definition}
 The length of such a sequence is bounded by the number of distinct boundary maps, an exponential function of $|\phi|$. And verifying whether the defining conditions of simple boundary maps, joins, limits, and shuffles hold is straightforward. Hence the set of fabricated boundary maps is decidable.  We say that $\phi$ \emph{occurs} in a fabricated boundary map $\partial$ if $\phi$ belongs to an MCS  or cluster in one of the $\partial$-labels, or inductively $\phi$ occurs in one of the previous boundary maps in the sequence for $\phi$ used to fabricate $\partial$.  

\begin{proposition} \label{thm:boundary} Let $\phi$ be a temporal formula.
The following are equivalent.
\begin{enumerate}
\item \label{eq:2} $\phi$ occurs in an open, fabricated boundary map $\partial$.
\item \label{eq:1} $\phi$ is satisfiable in some $\M$-model.
\end{enumerate}
\end{proposition}
An open boundary map is one summarising an open rectangle (of which $\mathbb R^2$ is a generic example), that is, a boundary map with domain $\{+, -\}$. Note that an open boundary map necessarily has no defects.

 Once we establish the proposition, it follows that  the validity of temporal formulas over $\M$ is decidable. In fact, \eqref{eq:2} can be decided nondeterministically with polynomial space, by searching recursively for a decomposition of $\partial$ as a join, limit, or shuffle of fabricated boundary maps, using tail recursion where possible \cite[Lemma 5.1]{HR18}.  The procedure to do this is shown in \Cref{alg:fabricated}.

\begin{algorithm}[H]
\caption{Nondeterministic procedure to decide whether $\partial$ is fabricated}\label{alg:fabricated}
\begin{algorithmic}
\Procedure{fabricated}{$\partial$}
{\bf choose} either

\State {\bf option} 0
\State \hspace*{\algorithmicindent}{\bf check} $\partial$ is simple

\State {\bf option} 1
\State \hspace*{\algorithmicindent}{\bf choose} $\partial_1, \partial_2$; {\bf check} they are boundaries

\State \hspace*{\algorithmicindent}{\bf check} their join ({\bf choose} some direction) is $\partial$; {\bf release} $\partial$

\State \hspace*{\algorithmicindent}{\bf check}  \Call {fabricated}{$\partial_1$}; {\bf tail-call}  \Call {fabricated}{$\partial_2$}

\State {\bf option} 2
\State \hspace*{\algorithmicindent}{\bf choose} $\partial_1, \partial_2, \partial_3, \partial_4$; {\bf check} they are boundaries

\State \hspace*{\algorithmicindent}{\bf check} $\partial$ is the limit ({\bf choose} direction) of $\partial_1$ using $\partial_2, \partial_3, \partial_4$; {\bf release}~$\partial$

\State \hspace*{\algorithmicindent}{\bf check}  \Call {fabricated}{$\partial_2}$,  \Call {fabricated}{$\partial_3$}; {\bf release} $\partial_2, \partial_3$  

\State \hspace*{\algorithmicindent}{\bf check} \Call {fabricated}{$\partial_4$}; {\bf tail-call}  \Call {fabricated}{$\partial_1$}

\State {\bf option} 3
\State \hspace*{\algorithmicindent}{\bf choose} $k \in \{1, \dots, |\phi|\}$, $\partial_1, \dots, \partial_k$ 

\State \hspace*{\algorithmicindent}{\bf check} they are boundaries including at least one one-point boundary

\State \hspace*{\algorithmicindent}{\bf check} $\partial$ is the shuffle of $\partial_1, \dots, \partial_k$; {\bf release} $\partial$

\State\hspace*{\algorithmicindent}{\bf for} $i=1, \dots, k$ {\bf do}

\State \hspace*{\algorithmicindent}\hspace*{\algorithmicindent}{\bf check}  \Call {fabricated}{$\partial_i$}
\State\hspace*{\algorithmicindent}{\bf end for}
\EndProcedure
\end{algorithmic}
\end{algorithm}

Hence \eqref{eq:2} is in NPSPACE, which by Savitch's theorem \cite{SAVITCH1970177} is equal to PSPACE. Thus validity of temporal formulas over $\M$ is in PSPACE. As noted already, the analogous result for $(\reals, \leq)\otimes(\reals, \leq)$ follows immediately, via the reduction given by applying $\P\psi \mapsto \psi \vee \P\psi$ and $\F\psi \mapsto \psi \vee \F\psi$ recursively to subformulas, and this validity problem is PSPACE-hard. Hence we obtain the following theorem.

\begin{theorem}[{\cite[Theorem 5.2]{HR18}}]\label{thm:leq}
The temporal logic of the frame $(\reals, \leq\nobreak)\otimes(\reals, \leq)$ is PSPACE-complete.   
\end{theorem}

To prove \Cref{thm:boundary}, we prove the more general equivalence between ($1'$)~$\partial$~being a fabricated boundary map and ($2'$)~existence of a rectangle model $h$ whose finite description is $\partial$.

The proof of ($1'$) $\Rightarrow$ ($2'$) is straightforward. We use dense valuations to construct rectangle models for simple boundary maps, and \Cref{fig:combine} shows how to synthesise rectangle models for joins, limits, and shuffles of more primitive boundary maps.   To be a bit more specific, a rectangle model for a simple boundary map may be constructed by using a dense valuation of clusters, that is, if a region (either the interior, or a segment of the boundary) is to be labelled by a cluster $c$ and $m\in c$, then at any point in the region there must be points in the region above and below labelled by $m$.   Joins and limits are straightforward (given that models may be variably dilated---in other words, reparameterised---horizontally/vertically, ensuring edges match up exactly).  Given rectangle models $h_1, \ldots, h_{k-1}$  for closed boundary maps $\partial_1, \ldots, \partial_{k-1}$ and a one-point boundary map $m$, an open rectangle model for a shuffle\label{p:shuffle} $\partial$ of the $\partial_i$ and $m$ may be constructed over the base $R=(0, 1)\times(0, 1)$ in steps.  Initially there is a single open `gap' $((0,1),(1, 0))$ in the line $x+y=1$, and no points are labelled.  At each stage, using a fair schedule, a gap is chosen and an $i \leq k$.  A copy of $h_i$ is used to label the rectangle whose diagonal is the closed, central third of the gap. (If $\partial_i$ is a one-point boundary map, then only the single point on the diagonal in the exact centre of the gap is labelled.)   The first and last third become new open gaps.  This is repeated $\omega$ times.  Finally, the Cantor set of unlabelled points on the line $x+y=1$ are all labelled $m$, all remaining unlabelled points where $x+y<1$ are given a label from $\partial(-)$ densely, and unlabelled points where $x+y>1$ are given labels from $\partial(+)$ densely, completing the definition of  the required rectangle model.

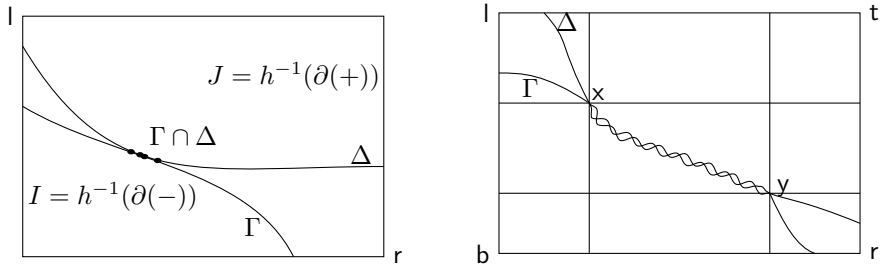
\begin{figure}
\begin{center}
\begin{tikzpicture}[xscale=.6,yscale=.4]

\draw (0,0) rectangle (8,8);

\draw  (6,0) .. controls (5,3)  and (2,3.2)   .. (0,5);

\draw (0,7) .. controls (2,2) and (4,3) .. (8,3);

\node at (2,2) {$I=h^{-1}(\partial(-))$};

\node at (6,6) {$J=h^{-1}(\partial(+))$};

\node at (5.1,1) {$\Gamma$};

\node at (7.5, 3.4) {$\Delta$};

\node [left] at (0,8) {$\ll$};
\node[right] at (8,0) {$\rr$};

\draw[fill] (2.7,3.33) circle [radius=2pt];
\draw[fill] (2.6,3.39) circle [radius=2pt];
\draw[fill] (2.4,3.49) circle [radius=2pt];
\draw[fill] (2.99,3.2) circle [radius=2pt];
\node[above right] at (2.6, 3.4) {$\Gamma\cap\Delta$};

\end{tikzpicture}
\hspace{.2in}
\begin{tikzpicture}[xscale=.6,yscale=.4,decoration={snake, amplitude=1.5, segment length=10}]
\draw (0,5)--(8,5);
\draw(0,2)--(8,2);
\draw (2,0)--(2, 8);
\draw(6,0)--(6, 8);
\draw (0,0) rectangle (8,8);
\draw(0,6) .. controls (.5, 6) and (1, 6) .. (2,5);
\draw(2, 5) .. controls (1.3, 7) and (1.6, 7)   ..(1,8);
\draw (6,2) .. controls  (6.3, 1) and (6.6, .2)   .. (7,0);
\draw (6, 2) .. controls  (6.3, 1.8) and (6.7, 1.8) .. (8, 1);
\node [left] at (0,8) {$\ll$};
\node[right] at (8,0) {$\rr$};
\node[left] at (0, 0) {$\bb$};
\node[right] at (8,8) {$\tt$};

\node at (2.2, 5.3) {$\x$};

\node at (6.3, 2.2) {$\y$};

\draw[decorate] (2, 5) .. controls (2.5, 3.5) and (3.5, 4.2) .. (6,2);
\draw[decorate] (6,2) .. controls (3.5, 3.5) and (2.5, 3.8) .. (2,5);

\node at (.7,5.5) {$\Gamma$};
\node at (1.5, 7.7) {$\Delta$};
\end{tikzpicture}

\caption{\label{fig:gamma} Maximal and minimal clusters with boundaries $\Delta,\Gamma$}
\end{center}
\end{figure}

The converse implication is  more intricate and only briefly outlined here (see \cite[Lemma 4.1]{HR18} for more details).   Let $h$ be a rectangle model. By reparameterising the coordinate axes, we may assume the rectangular domain of $h$ is bounded.  The proof that the rounded boundary map $\partial_h$ defined by $h$ is fabricated is by induction over the \emph{height} of the rectangle---the maximum length of an ordered chain of distinct clusters from $\partial_h(-)$ to $\partial_h(+)$.\footnote{The proofs appearing in \cite{HR18} and \cite{HM18} address only rectangular boundary maps, but the extension to rounded boundary maps is trivial.}  When the height is zero, $\partial_h(-)=\partial_h(+)$, so 
$\partial_h$ is simple, hence fabricated.    When the height is positive, we consider $I=h^{-1}(\partial_h(-))$, the set of points in the rectangle whose truth set belongs to $\partial_h(-)$---see the first diagram in \Cref{fig:gamma}.

 Modulo points on the rectangle edges, $I$ is a downward-closed set. The `upper boundary', $\Gamma$, of $I$, is topologically equivalent (in the usual, open-ball topology on $\reals^2$) to a closed line segment through the rectangle \cite[Lemma 2.11]{HR18}. Similarly, let $\Delta$ be the `lower boundary' of $J=h^{-1}(\partial_h(+))$.   If $\Gamma$ and $\Delta$ are disjoint then, as they are closed,  there is $\epsilon>0$ such that each point in $\Gamma$ is at least $\epsilon$ from each point of $\Delta$. So  $\partial_h$ is a join of the boundary maps of finitely many rectangle models of side at most $\frac\epsilon{\sqrt 2}$, of strictly smaller height, inductively fabricated; hence $\partial_h$ is fabricated.   Now suppose $\Gamma\cap\Delta$ is non-empty.    The set $\Gamma\cap\Delta$ is closed, so consider the  points $\x, \y$,  where $\x$ is the nearest point in $\Gamma\cap\Delta$ to $\ll$, and $\y$ is the nearest point in $\Gamma\cap\Delta$ to $\rr$---see the second diagram in \Cref{fig:gamma}.

  The rectangles $[\bb,\x],[\x,\tt], [\bb,\y], [\y,\tt]$ determine boundary maps of strictly smaller height, so they are fabricated.  
Consider the rectangle $[\ll\wedge\x,\ll\vee\x]$ with opposite corners $\ll$ and  $\x$.  There are only finitely many boundary maps, so let $\x_i$ be a sequence of points in this rectangle, strictly between $\Gamma$ and $\Delta$, converging to $\x$ such that the boundary map of the rectangle $[\ll\wedge\x_i, \ll\vee\x_i]$ is constant---call this boundary map $\partial_0$.   By the previous case (where $\Gamma$ does not meet $\Delta$) the map $\partial_0$ is fabricated. The boundary map of $[\ll\wedge\x, \ll\vee\x]$   is a southeastern limit of $\partial_0$, hence it is fabricated.   Similarly, the boundary map of $[\y\wedge\rr, \y\vee\rr]$ is fabricated.  

It remains to check that the boundary map of the rectangle $[\x\wedge\y, \x\vee\y]$ is fabricated.    If $\x=\y$ this is trivial, so assume not.  Let $\approx$ be the smallest equivalence relation over $\Gamma\cap\Delta$ including the successor relation (i.e. if $u, v\in\Gamma\cap\Delta$ and there is no $z\in\Gamma\cap\Delta$ strictly between $u$ and $v$ then $u\approx v$), including all pairs of points that differ in only one coordinate, and whose equivalence classes form closed sets (so if $\x_i\approx \y$ (all $i$) and $\x_i$ converges to $\x$ then $\x\approx \y$).   The closed and bounded set $\Gamma\cap\Delta$ is partitioned by these equivalence classes, so by topological properties, either there is only one equivalence class, or uncountably  many equivalence classes are singletons.    Each $\approx$-equivalence class $e$  also has a first and a last point $f(e)$ and $l(e)$, respectively. The rectangle $[f(e)\wedge l(e), f(e)\vee l(e)]$ defines a boundary map $\partial_e$ (and in some cases this must be a one-point boundary map) that can be shown to be fabricated by considering joins and limits of  boundary maps that, by previous cases, are fabricated. Hence each $\partial_e$ is also fabricated. (The argument here deviates from the proof of \cite[Lemma 4.1]{HR18}, and is more similar to the proof of \cite[Lemma 5.3]{HM18}.)  If there is only one $\approx$-equivalence class, we are done. Otherwise, the proof of \Cref{thm:boundary} is completed by showing that the boundary map of $[\x\wedge\y,\x\vee\y]$ is a shuffle of the boundary maps $\partial_e$.

\section*{Slower-than-light signals}
Now we consider the temporal logic of the frame $(\reals, <)\otimes(\reals, <)$, where the irreflexive accessibility relation is `can send a signal at strictly less than the speed of light'.  We can no longer give a direct description of the MCSs holding along the four light-lines bounding a rectangle, since distinct points on a light-line are now unordered.  What we can do, nevertheless, is record, at each point on the light-line, the cluster of MCSs that hold arbitrarily soon in the future (the \emph{upper cluster} of a point) and the cluster of MCSs that hold arbitrarily recently in the past (the \emph{lower cluster} of a point).  The function from points to  upper cluster (or to lower cluster) is monotone with respect to the parametric ordering of the light-line.  

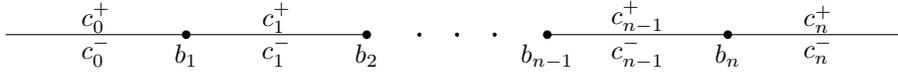
\begin{figure}[H]
\centering
\begin{tikzpicture}
\foreach \x in {0,1}
{
\draw[fill] (\tikzconst*\x+\tikzconst,0) circle [radius=0.05];
\node [above] at (\tikzconst*\x+.5*\tikzconst, -.1) {$c^+_\x$};
\node [below] at (\tikzconst*\x+.5*\tikzconst, 0.1) {$c^-_\x$};
\draw (\tikzconst*\x,0)--(\tikzconst*\x +\tikzconst,0);
}
\node [below] at (\tikzconst, 0) {$b_1$};
\node [below] at (\tikzconst*2, 0) {$b_2$};
\foreach \x in {3,4}
{
\draw[fill] (\tikzconst*\x,0) circle [radius=0.05];
\draw (\tikzconst*\x,0)--(\tikzconst*\x +\tikzconst,0);
}
\node [below] at (\tikzconst+\tikzconst*2, 0) {$b_{n-1}$};
\node [above] at (\tikzconst*3.5, -.1) {$c^+_{n-1}$};
\node [below] at (\tikzconst*3.5, 0.1) {$c^-_{n-1}$};
\node [below] at (\tikzconst+\tikzconst*3, 0) {$b_{n}$};
\node [above] at (\tikzconst*4.5, -.1) {$c^+_n$};
\node [below] at (\tikzconst*4.5, 0.1) {$c^-_n$};
\draw[fill] (2.5*\tikzconst-.5,0) circle [radius=0.02];
\draw[fill] (2.5*\tikzconst,0) circle [radius=0.02];
\draw[fill] (2.5*\tikzconst+.5,0) circle [radius=0.02];
\end{tikzpicture}
\caption{A bi-trace}\label{fig:bi-trace}
\end{figure}

\begin{definition}[Bi-trace]\label{def:bi-trace} A \emph{bi-trace} consists of two sequences of clusters $c_0^+\lesssim\ldots\lesssim c_n^+$ and $c_0^-\lesssim\ldots\lesssim c_n^-$, and one sequence of MCSs $b_1, \ldots, b_n$ (for some $n$) such that for all $i$:
\begin{enumerate}[(i)]
\item
 $c_i^-\lesssim c_i^+$ and $(c_i^-, c_i^+)\neq (c_{i+1}^-, c_{i+1}^+)$,
 \item
 $c_i^-\lesssim b_{i+1}\leq c_{i+1}^+$,  
\item\label{cond:3} there exists $m_i\in \mathsf{MCS}$ such that: $c_i^-\lesssim m_i\lesssim c_{i}^+$, \/all future defects of $m_i$ are passed up to $c_i^+$, and all past defects passed down to $c_i^-$ (such an $m_i$ is an \emph{interpolant} of $c_i^-$ and $c_{i}^+$), 
\item all future defects of $b_{i+1}$ are passed up to $c_{i+1}^+$ and all past defects of $b_{i+1}$ are passed down to $c_i^-$.
\end{enumerate} See \Cref{fig:bi-trace}.
\end{definition}

Intuitively, an interpolant of $c_i^-$ and $c_{i}^+$ is an MCS that can appear \emph{on} the light-line described by the bi-trace, on the segment between $c_i^-$ and $c_{i}^+$. %Condition \eqref{cond:3} of \Cref{def:bi-trace} ensures there is \emph{a priori} at least one such MCS available.

A \emph{bi-boundary} is like a boundary map, but uses bi-traces instead of traces.  
\begin{definition}[Bi-boundary]\label{def:bi-boundary} A \emph{bi-boundary} $\partial$ is a map from a rounded subset of    $\set{-, +, N, S, E, W, \bb, \ll, \rr, \tt}$ to clusters (for $-, +$), bi-traces (for $N, S, E, W$), and MCSs (for $\bb, \ll, \rr, \tt$) such that
\begin{itemize}
\item temporal ordering is respected, i.e. $\partial(+)$,  the final lower cluster of $\partial(N)$, and the final lower cluster of $\partial(E)$ are all equal (when defined), and $\partial(+)\lesssim \partial(\tt)$ (when $\partial(\tt)$ is defined),   plus a dual property for $\partial(-)$,
\item future defects of $\partial(+)$ are passed up to either an interpolant of the final lower and upper  clusters of $\partial(N)$ (provided $\partial(N)$ is defined), an interpolant of the final  lower and upper clusters of $\partial(E)$ (if defined) or to $\tt$ (if defined), plus a dual condition for past defects of $\partial(-)$.
\end{itemize} 
\end{definition}

 The definitions of joins, limits, shuffles, fabricated etc. above may be altered for bi-boundaries with no significant changes, hence the following theorem.
\begin{theorem}[{\cite[Theorem~6.1]{HM18}}]\label{thm:slower} 
The temporal logic of the frame $(\reals, <\nobreak)\otimes(\reals, <)$ is PSPACE-complete.
\end{theorem}

\section*{Real intervals}
There is a strong connection between these two-dimensional frames and interval logic, noted in \cite{SS02}.  Interval logics are based on frames consisting of strictly ordered pairs from a linear order.  The full interval logic has modalities for all thirteen of Allen's interval relations (after, meets, overlaps, ends, during, starts, and their converses, and equals), the full logic is undecidable and remains undecidable if we consider the logic of intervals $(x, y)$ with \emph{real} endpoints $x<y$ \cite{HaSh86}.  A large body of research has investigated the decidability and complexity of the interval logic where the set of modalities is restricted \cite{Bres19},  and over various linear orders. For example the modal interval logic with the two modalities `contained in' (i.e. during, starts, ends, or equals) and its converse, over real intervals, is known to be PSPACE-complete \cite{MPHS10}.    A  problem that appears to remain open is the decidability and complexity of the logic of intervals with real endpoints with the two modalities, strict during and its converse (see the open problem 3 of \cite{HaSh86}).   For comparison, the logic of intervals over \emph{discrete} time with strict during and its converse, is undecidable.
With some modification, we will use the argument above to prove decidability of the interval logic with strict during and its converse as modalities, over the reals. With further analysis the logic can be shown to be PSPACE-complete.

If we restrict the frame $(\reals, <)\otimes(\reals, <)$ to $\set{(x, y): x<y\in\reals}$ and treat its points as intervals, the accessibility relation is `each endpoint is strictly earlier', equivalently,  overlaps, meets, or before.   Instead, we restrict $(\reals, <)\otimes(\reals, <)$ to  $\set{(x, y): x+y>0}$  to obtain the open,  upper triangular frame $\c T$.  A point $(x, y)$ in this subframe represents the interval $(-x, y)$, more precisely, the map $(x, y)\mapsto (-x, y)$ is a frame isomorphism from $\c T$ to the frame of intervals $(x, y)$ with real endpoints $x<y$ under the during relation, i.e. $(x, y)\mathrel{\mathtt{dur}}(x', y')\iff x'<x<y<y'$.  
The temporal logic of this frame differs from the temporal logic of $(\reals, <)\otimes(\reals, <)$. In particular the past confluence axiom $\P\H p\rightarrow \H\P p$ is valid over $(\reals, <)\otimes(\reals, <)$ but not valid in $\c T$.

\begin{figure}
\begin{center}
\begin{tikzpicture}[scale=.35];

\draw[thick,->] (-1,0) -- (4,0);
\draw[thick,->] (0,-1) -- (0,4);
\draw[dashed](-3,3)--(3,-3);

\node at (3,3){$\bullet$};
\node [right] at  (3, 3){$(x, y)$};

%\node at (4,4){$\bullet$};
%\node [right] at  (4,4){$(x', y')$};

\draw (3,3)--(3,-3);
\draw(3,3)--(-3,3);

\node at (1.5,1.5){$T(x, y)$};

\end{tikzpicture}
\hspace{.2in}
\begin{tikzpicture}[scale=.5,decoration={snake, amplitude=1.5, segment length=10}]

\draw (0,4)--node [midway, above] {$\tau(N)$} (4,4);
\draw(4,4)-- node [midway, right]{$\tau(E)$} (4,0);
\draw[dashed](0,4)--(4,0);
\node at (4,4){$\bullet$};
\node[above] at (4, 4){$\tau(\tt)$};

\node at (3.3,3.5) {$\tau(+)$};
\draw[dashed, decorate] (2.4,4) arc (180:270:1.6cm);

%\draw[dashed](0,4)--(0, 5.5);

%\node at (.7,4.7){$\tau(\ll)$};
 % \draw[dashed] (1.5,4) arc (0:90:1.5cm) ;

%\draw[dashed] (4,0)--(5.5,0);

%\node at (4.5,0.5) {$\tau(\rr)$};
 % \draw[dashed] (5.5,0) arc (0:90:1.5cm);
\end{tikzpicture}
\vspace{-.15in}
\end{center}
\caption{\label{fig:tri} The  closed triangle frame  $T[x, y]$ and a closed triangle bi-boundary}
\end{figure}
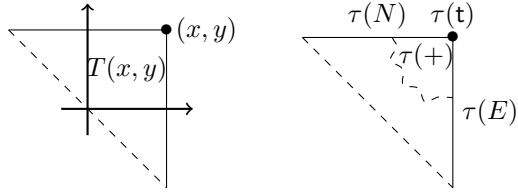

Let $x+y>0$.    The \emph{upper-closed triangle frame} $T[x, y]$ has worlds \[\set{(x', y'):  x'+y'>0,\; x'\leq x, \; y'\leq y}\] and accessibility inherited from $\c T$; see the first diagram in \Cref{fig:tri}. An \emph{open} (or \emph{semi-open}) triangle frame is obtained by requiring both (either) constraints $x'\leq x,\; y'\leq y$ to be strict.   Note that points on the line $x'+y'=0$ are excluded from the domain of triangle frames.  If $h$ is a propositional valuation over $T[x, y]$ and $x'<x,\;y'<y$ then $h$ induces not only a valuation  over $T[x', y']$ but also determines the upper cluster holding above but arbitrarily close to   points along the northern and eastern perimeter of $T[x', y']$.  
\begin{definition}[Closed triangle bi-boundary]\label{def:7}
A \emph{closed triangle bi-bound\-ary} $\tau$ is a map from $\set{+, N, E, \tt}$ to a cluster (for $+$), bi-traces (for $N, E$), and an MCS  (for $\tt$) satisfying the following consistency constraints. (See the second diagram in \Cref{fig:tri}.)
\begin{enumerate}
\item
The cluster $\tau(+)$ equals the final lower cluster of $\tau(N)$ and of $\tau(E)$, which 
is below $\tau(\tt)$.
  
\item\label{eq:F}
Future defects of $\tau(+)$ are passed
 up either to $\tau(\tt)$, or to an interpolant of $\tau(+)$   and the final upper cluster of $\tau(N)$, or to an interpolant of $\tau(+)$ and the final upper cluster of  $\tau(E)$.
 
\item
Past defects of $\tau(\tt)$ are passed down to $\tau(+)$.
\end{enumerate} 
\end{definition}

 Open and semi-open triangle
  bi-boundaries are similar. When $\tau$ is open, observe that none of the alternatives in \eqref{eq:F} are possible, so $\tau(+)$ has no future defects.  If $h$ is a propositional valuation over a triangle frame and $(x, y)$ is in the interior of the frame, we write $\tau_h[x, y]$ for the closed \tbb\ induced by the triangle $T[x, y]$. The semi-open  and open versions are denoted $\tau_h[x, y)$, $\tau_h(x, y]$, and $\tau_h(x, y)$.

 A \emph{simple} \tbb\ is one where $\tau(+)$ has no past defects,  and if the bi-traces $\tau(N)$ and/or $\tau(E)$ are defined then their lower clusters all equal $\tau(+)$.
 
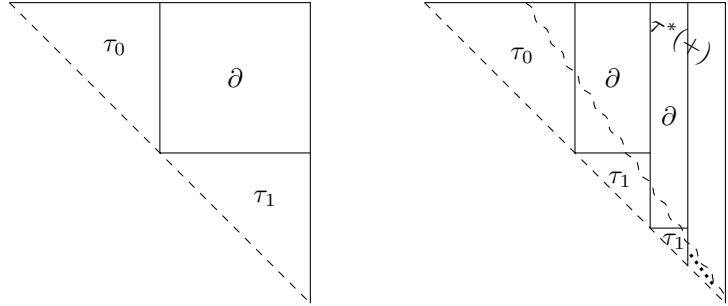
\begin{figure}[H]
\begin{center}
\begin{tikzpicture}[scale=1]
\draw[dashed](0,4)--(4,0);
\draw (0,4)--(4,4)--(4,0);
\draw (2,4)--(2,2)--(4,2);
\node at (1.4,3.4){$\tau_0$};
\node at (3,3){$\partial$};
\node at (3.4, 1.4){$\tau_1$};

\end{tikzpicture}
\hspace{.5in}
\begin{tikzpicture}[scale=1, decoration={snake, amplitude=.8}]]
\draw[dashed](0,4)--(4,0);
\draw (0,4)--(4,4)--(4,0);
\draw (2,4)--(2,2)--(3,2)--(3,4);
\draw (3,2)--(3,1);
\draw (3.5, 4)--(3.5, .5);
\draw (3,1)--(3.5,1);
\node at (1.3,3.3){$\tau_0$};
\node at (2.5,3){$\partial$};
\node at (3.25,2.5){$\partial$};
\node at (2.6, 1.65){$\tau_1$};
\node at (3.32, .82){$\tau_1$};
\node[rotate = -30] at (3.4, 3.5) {$\tau^*(+)$};

\draw[dashed, decorate](1.35,4)--(4.005,0.09);

\draw[dotted, very thick] (3.535,.65)--(3.85,0.23);

\end{tikzpicture}
\vspace{-.15in}
\end{center}
\caption{\label{fig:tri2} The join $\tau_0\oplus\partial\oplus\tau_1$ and a southeastern limit $\tau^*$ of $\tau_0$}
\end{figure}

Given two closed triangle bi-boundaries  $\tau_0, \tau_1$, and a  bi-boundary $\partial$ with domain $\set{-, +,  N, S, E, W, \ll,\rr,\tt}$, if $\tau_0(E)=\partial(W),\;\tau_0(\tt)=\partial(\ll),\; \partial(S)=\tau_1(N),\allowbreak\partial(\rr)=\tau_1(\tt)$, then we may form the closed triangle \emph{join} $\tau_0\oplus\partial\oplus\tau_1$ as illustrated in the first diagram of \Cref{fig:tri2}.    Observe that the domain of $\partial$ omits $\bb$, so it is not rectangular, but it \emph{is} rounded. Triangle joins that are open or semi open (at the north and/or east) are defined similarly.

The \tbb\  $\tau^*$ is a \emph{southeastern limit} of $\tau_0$ if
\begin{itemize}
\item  $\tau_0=\tau_0\oplus\partial\oplus\tau_1$ for some $\partial, \tau_1$,
\item
the lower cluster of the bi-trace $\partial(E)$ constantly $\tau^*(+)$,
\item
$\tau_0$ agrees with $\tau^*$ over $\set{+, N}$ (allowing both to be undefined on $N$), %\/ $\tau^*$ is defined on $E$, \nb{surely there are E-open limits too}
\item
$\tau^*(E)$, if defined, is a bi-trace where every lower cluster equals $\tau_0(+)$ (hence $\tau^*(+)= \tau_0(+)=\tau_1(+)$).
\end{itemize} See the second diagram in \Cref{fig:tri2}, where the wavy line indicates the implied boundary of the region labelled with $\tau^*(+)$ %and every future defect of $\tau_0(+)$ is passed up to either $\tau^*(N)$ (when defined) or $\tau^*(E)$.\nb{this part is already in the definition of bi-boundaries}
 \emph{Northwestern} limits are similar. 

 If $\tau_1,\ldots, \tau_{k}$ are closed triangle bi-boundaries then $\tau$ is a \emph{shuffle} of them if every past defect of $\tau(+)$ is passed down to $\tau_i(\tt)$ (some $i \leq k$) and every future defect of $\tau_i(\tt)$ is passed up to $\tau(+)$ (all $i\leq k$). Observe that a one-point \tbb\ is not required for a shuffle since the diagonal boundary edge is excluded from the domain of the triangle, in contrast to the situation for shuffles of rounded bi-boundaries.

 The set of \emph{fabricated} triangle bi-boundaries is the closure of the set of simple triangle bi-boundaries under joins (using fabricated rounded bi-boundaries), limits, and shuffles.

 \begin{proposition}\label{thm:4}
 Let $\tau$ be  a closed (respectively open, semi-open north, semi-open east)  triangle bi-boundary map.  The following are equivalent.
 \begin{enumerate}
 \item $\tau$ is fabricated.
 \item There is a valuation $h$  over a triangle frame $T$ with $(1, 1)$ in the interior of $T$ such that $\tau=\tau_h[1,1]$ (respectively, $\tau=\tau_h(1,1), \tau_h(1,1], \tau_h[1,1)$).
 \end{enumerate}
 \end{proposition}

As with rectangle models and boundary maps, it is simple to check that  any fabricated triangle bi-boundary may be obtained as $\tau_h[x, y]$, $\tau_h[x, y)$, $\tau_h(x, y]$, or $\tau_h(x, y)$, for some valuation $h$. (To be completely accurate, for any fabricated triangle bi-boundary $\tau$ there is an $h$ such that the induced $\tau_h$ agrees with $\tau$ everywhere \emph{except possibly} the upper clusters of $\tau(N)$/$\tau(E)$. Importantly, this technicality disappears in the case of open triangles---those corresponding to $\c T$.)

Conversely, for any valuation $h$ over some triangle frame $T$ where $(1, 1)$ is in the interior, we must show that $\tau_h[1, 1]$, $\tau_h[1, 1)$, $\tau_h(1, 1]$, and 
$\tau_h(1, 1)$ are fabricated. We focus on $\tau_h[1,1]$ below; the other cases are similar.  The proof is by induction over the \emph{size} of the triangle model---the number of distinct clusters witnessed in the open triangle $T(1,1)$.  Let $\lambda$ be the line $\set{(x, y):x+y=0}$.  For $(x, -x)\in\lambda$ and $x<y,\; -x<z$ we write $_*[(x,-x), (y, z)]$ for the `rounded rectangle' $[(x, -x), (y, z)]\setminus\set{(x, -x)}$.  By the bi-boundary equivalent of \Cref{thm:boundary}, the bi-boundary $\partial_h(_*[(x, -x), (y, z)])$ is fabricated.

 For the base case, the single cluster  $\tau_h(+)$ covers the entirety of $T(1,1)$. The cluster $\tau_h(+)$ must have no past defects, and $\tau_h(N)$ and $\tau_h(E)$, if defined, must have all lower clusters equal to $\tau_h(+)$. Future defects of $\tau_h(+)$ must be passed up. Hence $\tau_h$ is simple, so fabricated.
 
   Now suppose that more than one cluster is witnessed in the interior, so  $J=h^{-1}(\tau_h(+)) \cap T(1,1)$ is a proper, upward-closed subset of the interior of the triangle.  Let $\Delta$ be the `lower boundary' of $J$. Formally, $\Delta$ is given by first taking the boundary of $J$ viewed as a subspace of $T(1,1)$, then taking the closure of that within $\reals^2$.
   For topological reasons, $\Delta$ is topologically equivalent to a line segment, closed by definition.

 If $\Delta$ does not meet $\lambda$ then since the sets are closed, and $\Delta$ is bounded, there is $\epsilon>0$ such that all points in $\Delta$ are at least $\epsilon$ from any point in $\lambda$.  The triangle $T[1,1]$ (whose orthogonal sides are of length 2) may be divided into at most $t=\lceil \frac {2\sqrt 2} \epsilon \rceil$ triangles $T(x, y)$, with $(x, y)$ below $\Delta$, of side at most $\frac{\epsilon}{\sqrt 2}$ and at most $t-1$ rectangles. Then since each triangle has strictly smaller size, by induction and by the result for bi-boundaries of rectangles, $\tau_h$ is fabricated.
 
\begin{figure}
\begin{center}
\begin{tikzpicture}[scale=.25,decoration={snake, amplitude=1.5, segment length=10}]

\draw (-10,10)--(10,10)--(10,-10);
\draw[dashed](-10,10)--(10,-10);

\node at (-3, 10){$\bullet$};
\node[above] at (-3,10) {$(x, 1)$};
%\node at (-3,3) {$\bullet$};
\node[below] at (-3, 3){$(x, -x)$};

%\node at (5, -5){$\bullet$};
\node[below] at (5, -5) {$(y, -y)$};

\node at (10, -5){$\bullet$};
\node[right] at (10, -5){$(1, -y)$};

\node at (5,3){$\bullet$};
\node[above] at (5, 3){$(y, -x)$};

\node at (8, 8) {$\tau(+)$};

%\node at (-3.6,8){$\tau(+)$};
%\node [fill=white] at (-4.5, 7){$\Delta$};
\node[right] at (10, 3){$(1, -x)$};

%\node at (8, -4.7){$\tau(+)$};

 \draw[decorate]  (-6,10,0) -- (-3,3) node [midway, fill=white] {{$\Delta$}};

%\draw[decorate(-6,10)--++(-3,3) node [midway,fill=white]{a};
\draw[decorate](5,-5)--(10, -6);
\draw(-3,10)--(-3,3)--(10,3);

\draw (10, -5)--(5, -5)--(5, 3);
\node[above] at (10, 10){$(1, 1)$};
\node[above] at (-10, 10){$(-1, 1)$};

\end{tikzpicture}
\end{center}
\caption{\label{fig:decomp}A decomposition of $\tau=\tau[1, 1]$ as $\tau_h[x, 1]\otimes\partial_h(_*[
(x, -x), (1, 1)] ) \otimes(\tau_h[y, -x]\otimes\partial_h (_*[ (y, -y), (1, -x)])\otimes\tau_h[1, -y])$. The wavy line indicates part of the lower boundary $\Delta$ of $\tau(+)$.}
\end{figure}
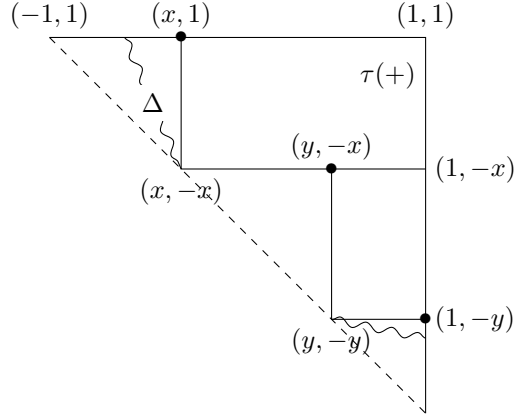

 So assume $\Delta\cap \lambda$ is non-empty.   Let $x\leq y$ be the infimum and supremum respectively of $\set{z:(z, -z)\in\Delta}$---see \Cref{fig:decomp}.  Since $\Delta\cap\lambda$ is closed, we know $(x, -x), \allowbreak(y, -y)\in\Delta\cap\lambda$.   Suppose $x\neq -1$. Then first we wish to fabricate $\tau_h[x, 1]$. We may assume there are no points of $\Delta$ directly north of $(x, -x)$ (otherwise we could fabricate $\tau_h[x, 1]$ using a triangle join and the $\Delta\cap\lambda=\emptyset$ cases, for triangles and rectangles). By this assumption, and as there are only finitely many triangle bi-boundaries, there is an increasing sequence $x_i$ converging to $x$ such that: %$(x_{i+1}, -x_i)$ is below $\Delta$, 
 the sequences $\tau_h[x_{i+1}, -x_i]$, $\tau_h[x_i, 1]$, and $\partial_*[(x_i, -x_i), (x_{i+1}, 1)]$ are all constant, and the bi-boundary $\partial_*[(x_i, -x_i), (x_{i+1}, 1)](E)$ is the single cluster $\tau_h(+)$.  %Since $\tau_h[x_{i+1}, -x_i]$ is of smaller size, inductively it is fabricated, and 
Since $\Delta$ is bounded away from $\lambda$ in $\tau_h[x_i, 1]$, by the $\Delta\cap\lambda=\emptyset$ cases   all three of $\tau_h[x_{i+1}, -x_i]$, $\tau_h[x_i, 1]$, and $\partial_*[(x_i, -x_i), (x_{i+1}, 1)]$ are fabricated.  Then $\tau_h[x, 1]$ is a southeastern limit of the constant \tbb\ $\tau_h[x_i, 1]$. Hence $\tau_h[x, 1]$ is fabricated.  %Also, $\tau_h[x_0, 1]$ is fabricated by the , hence $\tau_h[x,1]$ is fabricated.    Similarly,  $\tau_h[1, -y]$ is fabricated. 

 If $y=x$ then $\tau$ is the fabricated join $\tau_h[x, 1]\otimes  \partial_h ( _*[(x, -x), (1,1)])\otimes\tau_h[1, -x]$.    So assume $x<y$. Then $\tau$ is the join 
$\tau_h[x, 1]\otimes\partial_h(_*[
(x, -x), (1, 1)] ) \otimes(\tau_h[y, -x]\otimes\partial_h (_*[ (y, -y), (1, -x)])\otimes\tau_h[1, -y])$, as shown in \Cref{fig:decomp}.
So  it remains to show that   $\tau_h[y, -x]$ is fabricated.

Let $\vartriangleleft$ be the northwest-to-southeast ordering on $\Delta\cap\lambda$.  Let $\approx$ be the smallest equivalence relation over $\Delta\cap\lambda$ such that 
\begin{enumerate}[(a)]
\item
 $\approx$ includes the successor relation  over $(\Delta\cap \lambda, \vartriangleleft)$ (i.e. if  $(x', -x'), (y', -y') \in\Delta\cap \lambda$ and  there is no $z$ strictly between $x'$ and $y'$ such that $(z, -z) \in\Delta\cap\lambda$
   then $(x', -x')\approx (y', -y')$),
   \item each $\approx$-equivalence class is closed in $\reals^2$ (so if $y_i$ is a sequence converging to $y'$ and $(x, -x)\approx (y_i, -y_i)$ (all $i$) then $(x', -x')\approx (y', -y')$).
   \end{enumerate}
 By minimality of $\approx$, the equivalence classes are convex with respect to $\vartriangleleft$. It follows that, as with rounded boundary maps, there is either a single $\approx$-equivalence class or an uncountable number of them.      

 Each (non-singleton) $\approx$-equivalence class $e$ determines an upper-closed triangle $T[\bigvee e, -{\bigwedge e}]$ where $\bigvee e$ and $\bigwedge e$ are the supremum and infimum, respectively, of the $x$-coordinates of $e$. As with rounded boundary maps, by considering joins and limits of fabricated triangle bi-boundaries, the bi-boundary of each $T[\bigvee e, -{\bigwedge e}]$ can be shown to be fabricated.  If there is a single $\approx$-equivalence class, we are done. Otherwise, the proof of \Cref{thm:4} is completed by showing that  $\tau_h[y, -x]$ is a shuffle of triangle bi-boundaries $\tau_h[\bigvee e, -{\bigwedge e}]$ where $e$ ranges over all non-singleton $\approx$-equivalence classes.

\begin{theorem}
The temporal logic of the frame of intervals $\{(x, y) \in \reals^2\mid x<y\}$ under the \emph{during} relation, $(x, y)\mathrel{\mathtt{dur}}(x', y')\iff x'<x<y<y'$, is PSPACE-complete.
\end{theorem}
 PSPACE-hardness of the temporal logic follows from the PSPACE-hardness of the modal logic of $(\reals, <)\otimes(\reals, <)$ \cite{Shap05}.   By \Cref{thm:4}, we now know that the temporal formulas satisfiable in the \emph{during} frame of intervals is decidable---it suffices to  check if a given formula $\phi$ occurs in an open, fabricated \tbb.    To obtain a polynomial space upper bound on the complexity, note first that just like for boundary maps and bi-boundary maps, a \tbb\ can be stored in polynomial space, in terms of $|\phi|$. Then just as for \Cref{thm:leq} and \Cref{thm:slower}, whether a triangle bi-boundary map is fabricated can be decided nondeterministically with polynomial space by searching recursively for a decomposition as a join, limit, or shuffle, using tail recursion where possible.  The procedure to do this is shown in \Cref{alg:fabricated2}.

 \begin{algorithm}[!htb]
\caption{Nondeterministic procedure to decide whether $\tau$ is fabricated}\label{alg:fabricated2}
\begin{algorithmic}
\Procedure{fabricated}{$\tau$}
{\bf choose} either

\State {\bf option} 0
\State \hspace*{\algorithmicindent}{\bf check} $\tau$ is simple

\State {\bf option} 1
\State \hspace*{\algorithmicindent}{\bf choose} $\tau_0, \partial, \tau_1$
\State \hspace*{\algorithmicindent}{\bf check} $\tau_0, \tau_1$ are triangle bi-boundaries; {\bf check} $\partial$ is a bi-boundary

\State \hspace*{\algorithmicindent}{\bf check} $\tau =\tau_0\oplus\partial\oplus\tau_1$; {\bf release} $\tau$

\State \hspace*{\algorithmicindent}{\bf check}  \Call {fabricated}{$\tau_0$}, \Call {fabricated}{$\tau_1$}; {\bf tail-call}  \Call {fabricated}{$\partial$}

\State {\bf option} 2
\State \hspace*{\algorithmicindent}{\bf choose} $\tau_0, \partial, \tau_1$

\State \hspace*{\algorithmicindent}{\bf check} $\tau_0, \tau_1$ are triangle bi-boundaries; {\bf check} $\partial$ is a bi-boundary

\State \hspace*{\algorithmicindent}{\bf check} $\tau$ is the limit ({\bf choose} direction) of $\tau_0$; {\bf release}~$\tau$

\State \hspace*{\algorithmicindent}{\bf check}  \Call {fabricated}{$\tau_1}$; {\bf release} $\tau_1$;

\State \hspace*{\algorithmicindent}{\bf check} \Call {fabricated}{$\tau_0$}; {\bf tail-call}  \Call {fabricated}{$\partial$}

\State {\bf option} 3
\State \hspace*{\algorithmicindent}{\bf choose} $k \in \{1, \dots, |\phi|\}$, $\tau_1, \dots, \tau_k$; {\bf check} they are bi-boundaries

\State \hspace*{\algorithmicindent}{\bf check} $\tau$ is the shuffle of $\tau_1, \dots, \tau_k$; {\bf release} $\tau$

\State\hspace*{\algorithmicindent}{\bf for} $i=1, \dots, k$ {\bf do}

\State \hspace*{\algorithmicindent}\hspace*{\algorithmicindent}{\bf check}  \Call {fabricated}{$\tau_i$}
\State\hspace*{\algorithmicindent}{\bf end for}
\EndProcedure
\end{algorithmic}
\end{algorithm}

In \Cref{alg:fabricated2}, a call of the form $\textsc{fabricated}(\partial)$ indicates a call to the bi-boundary analogue of \Cref{alg:fabricated} needed to prove \Cref{thm:slower}, which is \cite[Algorithm~1]{HM18}.
 
\section*{Open problems}
\begin{enumerate}
\item  Find an axiomatisation of the valid temporal formulas over $(\reals, \leq)\otimes(\reals, \leq)$.
\item  Generalise the material above to find more general conditions where the temporal logic of  $(W, R)\otimes(W', R')$ is decidable.
\item  Determine the decidability of the temporal logic of Minkowski spacetime of higher dimensions. We conjecture this logic is undecidable for three or more spacetime dimensions.

\end{enumerate}

\newcommand{\etalchar}[1]{$^{#1}$}
 \def\www{/\allowbreak}

%\bibliographystyle{alpha}
%\bibliography{/Users/robin/Dropbox/tex/bibinputs/robin}
\end{document}